\documentclass[12pt]{article}
\usepackage{tabularx}
\usepackage[usenames,dvipsnames,table]{xcolor} 
\usepackage[justification=centering]{caption}
\usepackage{graphicx}
\usepackage{float}

\usepackage{bm}
\usepackage{tikz}

\def\opt{{\rm OPT}}
\def\eop{\hfill\fbox{}\medskip}

\newtheorem{Theorem}{Theorem}

\newtheorem{Property}{Property}

\long\def\comment#1\endcomment{}
\date{}

\begin{document}
\begin{center}
\LARGE
Scheduling a single machine with compressible jobs to minimize
maximum lateness

\vspace{7mm}
 
\normalsize

Nodari Vakhania  \\
Centro de Investigaci\'on en Ciencias, Universidad Aut\'onoma del Estado de Morelos, \\  Cuernavaca 62209, Mexico\\
email: nodari@uaem.mx  \\[3ex]
Frank Werner  \\
Faculty of Mathematics, Otto-von-Guericke University Magdeburg, \\  PSF 4120, 39106 Magdeburg, Germany  \\
email: frank.werner@ovgu.de  \\[3ex]
Alejandro Reynoso  \\
Centro de Investigaci\'on en Ciencias, Universidad Aut\'onoma del Estado de Morelos, \\  Cuernavaca 62209, Mexico\\
email: amlh26c7@gmail.com  \\[3ex]
\end{center}


\begin{abstract} {The problem of scheduling non-simultaneously 
released jobs with due dates on a single machine with the 
objective to minimize the maximum job lateness is known to be
strongly NP-hard. Here we consider an extended model in which 
the compression of the job processing times is allowed. The
compression is accomplished at the cost of involving additional 
emerging resources, whose use, however, yields some cost. With a given 
upper limit $U$ on the total allowable cost, one wishes to minimize 
the maximum job lateness. It is clear that, by using the available 
resources, some jobs may complete earlier and the objective function 
value may respectively be decreased. As we show here, for minimizing the 
maximum job lateness, by shortening the processing time of some specially 
determined jobs, the objective value can be decreased. Although the generalized 
problem is  harder than the generic non-compressible version, given a 
``sufficient amount'' of additional resources, we can solve the problem 
optimally. We determine the compression rate for some specific jobs
and develop an algorithm that obtains an optimal solution. Such an approach can be
beneficial in practice since the manufacturer can be provided with an information
about the required amount of additional resources in order to solve the 
problem optimally. In case the amount of the available additional resources is 
less than used in the above solution, i.e., it is not feasible, it is transformed
to a tight minimal feasible solution.}
\end{abstract}


\newpage
 
\section{Introduction}

The main part of the scheduling literature deals with deterministic problems, where all data including the processing times are fixed values. However, in many practical problems the processing time of a job can be controlled e.g. by the amount of an allocated resource used 
the processing time can be compressed from a standard value to a smaller value. In this paper, we consider such a problem with compressible processing times which can be described as follows. 
 
The $n$ jobs 
$1, 2, \ldots, n$ from a set $J$ are to be processed on a single machine.
Each job $j$ is available for processing from its {\em release time} $r_j$ 
and has the desired completion time or {\em due date} $d_j$. 
For each job $j$, we have an initially given {\em processing time} $a_j$ 
and the {\em cost} for the unitary compression $c_j$. If the processing time 
of job $j$ is compressed by $x_j(S)$  time units 
in the schedule $S$, then its real processing time is $$p_j(S)=a_j-x_j(S), $$ and 
the cost for this compression is $x_j(S) c_j$. There is an upper limit 
$U\ge 0$ on the total compression cost. Hence, the total cost in a feasible 
schedule $S$ has to be not larger than $U$, i.e., 
\begin{equation}\label{1}
\sum_{j}x_j(S) c_j \le U.
\end{equation} 
We shall refer to an $n$-component vector $(x_1(S),\dots,x_n(S))$ as the 
{\em compression vector} for a schedule $S$. The {\em completion time} of
 job $j$ in the schedule $S$ is $$f_j(S)=s_j(S)+p_j(S), $$ where $s_j(S)$ is the 
{\em starting time} of job $j$ in that schedule. The {\em lateness} of job $j$
in the schedule $S$ is $$L_j(S)= f_j(S)-d_j,$$ and the lateness $L(S)$ of schedule $S$ is
the maximum job lateness in it, i.e., $$L(S)=\max_j L_j(S). $$ The aim is to find an 
optimal feasible schedule, i.e., a feasible schedule with the minimum lateness $L^{\opt}$.

This setting is motivated by real-life scenarios, where additional resources are
available for processing the jobs. In case an additional resource is 
used for job $j$, the processing requirement of that job reduces accordingly. 
The cost for incorporating an additional resource for job $j$ is reflected by the 
parameter $c_j$. The total budget for the additional resources is limited 
by a constant $U$, which makes it harder to solve the problem (similarly, as 
bounding the knapsack capacity makes the KNAPSACK problem hard). Note that
the problem $1|r_j, compressible(U)|L_{\max}$ is a special case of the generic 
version $1|r_j|L_{\max}$ with $U=0$. Hence, the former problem is as hard as
the latter generic one which is known to be strongly NP-hard \cite{GJ}. 
The problem can also be seen as a bi-criteria optimization problem with two
contradictory objective criteria, to minimize the maximum lateness and to
minimize the total compression cost. It is not difficult to see that the
corresponding Pareto optimization problem remains strongly NP-hard. 
Among the feasible schedules with a given total cost $U$, finding one
with the minimum lateness is clearly as hard as problem $1|r_j|L_{\max}$ (in
particular, for $U=0$, we have an instance of problem $1|r_j|L_{\max}$). 
Likewise, among the feasible solutions with a given maximum job lateness, it 
is hard to find one with the minimum total compression cost. 
 
There are a number of exact exponential-time implicit enumeration  
algorithms for problem $1|r_j|L_{\max}$. However, it is not easy
to arrive at an efficient enumerative method for the extended problem
$1|r_j, compressible(U)|L_{\max}$ since it is not easy to enumerate 
all feasible ways for the compression of the job processing times;  here, 
for each possible selection of compressed processing times, an exact 
algorithm for problem $1|r_j|L_{\max}$ has to be invoked. Hence, it is 
unlikely that one can arrive at an easy and efficient exact solution method 
for the general setting $1|r_j, compressible(U)|L_{\max}$. Given this, 
a polynomial-time algorithm that gives optimal solutions under some
reasonable conditions and reasonable sub-optimal solutions would clearly be of
an interest. In this paper,  we describe such a polynomial-time method 
considering two alternative practical versions of the problem 
$1|r_j,compressible(U)|L_{\max}$.  

The method that we propose provides the manufacturer with an 
amount of additional resources; given this amount of additional resources, 
the method obtains an optimal solution to the problem. If the amount of 
additional resources is a priory fixed and the total compression 
cost of our solution is greater than $U$ (i.e., it is not feasible
for problem $1|r_j, compressible(U)|L_{\max}$), we transform it to a 
feasible solution which is optimal for the corresponding instance of the 
generic problem $1|r_j|L_{\max}$. Our transformation delivers a tight minimal 
feasible solution problem $1|r_j, compressible(U)|L_{\max}$, in the sense that 
by increasing the processing time of a job by some amount of time, the maximum 
job lateness in the resultant solution will increase by the same amount. In
particular, for a chosen compression vector 
$(x_1,\dots,x_n)$, an instance of problem $1|r_j|L_{\max}$ with the job 
processing times $a_1-x_1,\dots, a_n-x_n$ is naturally associated 
with the initially given instance of problem $1|r_j, 
compressible(U)|L_{\max}$. The two algorithms that we propose, for
a given instance $I$ of problem $1|r_j, compressible(U)|L_{\max}$, find a
compression vector with the associated instance $I'$ of problem $1|r_j|L_{\max}$
and a solution $S$ respecting that vector, which is optimal for the instance $I'$. 
If the solution $S$ satisfies condition (\ref{1}), then it is also optimal for 
instance $I$. Otherwise, it is transformed to a tight feasible solution, as 
indicated above. 

Our first algorithm is pseudo-polynomial, its running time depends linearly on
 the maximum job processing time $p_{\max}$.  The second one is polynomial and 
 runs in time $O(n^2\log n)$. It applies a similar strategy as the pseudo-polynomial
one but the search is organized in a different way so that the dependence
on  $p_{\max}$ is avoided.  
 
The remainder of the paper is organized as follows. In Section 2, 
we give a brief overview of the related literature. Section 3 presents 
some necessary preliminaries. Section 4 consists of three subsections.
In Section 4.1 we describe an algorithm that provides us with our initial solution. 
In Section 4.2, we discuss the proposed approach in details and give our 
pseudo-polynomial time algorithm. In Section 4.3 we describe our polynomial-time 
algorithm. Finally, Section 5 gives a few concluding remarks.

\section{Literature Review and Applications}

In this section we survey the earlier work that we found in the literature for 
scheduling problems using some kind of the control or compression of job processing times, where 
they are not constant. Since the 1980ies, a lot of papers dealing with such a scheduling environment appeared so that we mention here only a few of the most relevant ones. 
One of the first papers considering scheduling problems with so-called controllable processing times was given by Vickson \cite{vic} in 1980. He considered two single machine problems and presented a polynomial algorithm as well as a complexity result. While the problem of minimizing total completion time can be formulated as an assignment problem and thus be solved in $O(n^{2.5})$ time, the consideration of job weights makes the problem $NP$-hard. Such problems with controllable processing times have not only theoretical but also practical importance. Janiak \cite{jan} described a real-world problem in the context of steel production,  where batches of ingots have to be preheated before they are hot-rolled in a blooming mill. Both the preheating and rolling times are indirectly proportional to the gas flow intensity. Another application arises in a machine tooling environment when the processing time of a job  depends on the feed rate and the spindle speed for the particular operations. Using our model for such applications, the manufacturer can be provided 
by the required amounts of gas portions for an optimal production process. Alternatively,
if the amount of gas is insufficient, then our method would adopt the solution to 
a feasible minimal solution that uses precisely the available amount of resource (gas).

The paper \cite{hoo} considers single machine scheduling problems with controllable processing times and minimizing the maximum job cost. For this problem, several polynomial time results are derived. On the other hand, for the corresponding problem with minimizing the total weighted completion time of the jobs an NP-hardness result has been presented. In \cite{zdr}, the single machine problem with given release dates, the criterion of minimizing maximum lateness subject to linearly controllable processing times has been considered. For this problem, a polynomial time approximation scheme was derived with a worst case approximation ratio of 11/6. 

A very detailed excellent survey on scheduling problems with controllable processing times for the period up to 2007 has been given in \cite{sha}, where 113 papers are surveyed. Most reviewed results refer to single machine problems. In separate subsections, the authors present results for the single machine problem with minimizing the maximum penalty term, total weighted completion time of the jobs, the weighted number of tardy jobs, batch scheduling problems, and due date assignment problems. In addition, for multi-machine problems, detailed tables with polynomial algorithms, approximation algorithms and complexity results are given.

Shakhlevich et al. \cite{shak}  deal with bi-criteria single machine scheduling  with controllable processing times, where the maximum cost depending on the job completion times and  the total compression cost should be minimized. This problem is reduced to a series of linear programs defined over the intersection of a submodular polyhedron with a box, and then a greedy algorithm is applied which is faster than earlier ones. Later, the same authors \cite{shi} consider a single machine problem with controllable processing times as well as given release date and  deadlines for each job subject to the minimization of the total cost for reducing the processing times. They reformulate the problem as a submodular optimization problem and develop a recursive decomposition algorithm. 

In the paper \cite{yin}, a single machine scheduling problem with controllable processing times and learning effect is investigated. The objective is to minimize  a cost function, containing makespan, total completion (waiting) time, total absolute differences in the completion (waiting) times and total compression cost.   The resulting problem is formulated as  an assignment probem and can thus be solved in polynomial time.
Shabtay and Zofi \cite{sha2} consider a single machine scheduling problem with controllable processing times  with an unavailability period and the objective of minimizing the makespan. For the case when the processing times are convex decreasing functions of the amount of the allocated resource, they show that the problem is NP-hard and present both  a constant factor approximation algorithm as well as a fully polynomial time approximation scheme. In the paper \cite{kay}, the single machine problem with controllable processing times and no inserted idle times with minimizing total tardiness and earliness is considered. After presenting a mathematical model, several heuristic approaches are suggested and compared on instances of different sizes. Recently, Luo and Zhang \cite{luo} considered a single machine problem, where in addition to controllable processing times also the setup times can be controlled. For the case of job and position-dependent workloads and minimizing the makespan, the authors show that this problem can be optimally solved in polynomial time. 

 A special situation happens when processing times can be compressed.  Cheng et al. \cite{cheng} consider a single-machine scheduling problem with common due date assignment and compressible processing times. The authors consider two variants of due-date assignment methods with the goal to determine an optimal job sequence, the optimal due dates  and the optimal compressions of the processing times in order to minimize a total penalty function. Such an approach somewhat resembles ours, where we determine optimal job compression rates to minimize the maximum lateness.
A single machine problem with randomly compressible processing times, which may result e.g. from the introduction of a new technology, has been considered by Qi et al. \cite{qi}.
Considering the cost for the breakdown and the compressible processing times, it is proven that under certain conditions, the optimal sequence satisfies the V-shape property.

Cao et al. \cite{cao} deal with three scheduling problems, where jobs can be rejected or the processing times can be discretely compressed.  All three problems considered are NP-hard. For two of these problems pseudo-polynomial dynamic programming algorithms and FPTASs are presented, and for the third problem a greedy heuristic is given.
Zhang and Zhang \cite{zhang} consider a scheduling problem with identical parallel machines to minimize the makespan subject to a constraint on the total compression cost. They suggest a pseudo-polynomial dynamic programming algorithm and an FPTAS. Peng et al. \cite{peng} consider  a single resource scheduling problem with compressible processing times. The goal is to minimize  the length of the delay time  and the number of compressed tasks. For this problem, a heuristic algorithm is presented and tested.

A related class of scheduling problems are those with deterioration, where the processing time depends on the starting time of the job, i.e., the later a job starts, the larger is its processing time. In \cite{kel}, a single machine problem with minimizing the makespan  subject to a cumulative deterioration is considered, where the processing conditions can be restored by a single maintenance. The authors investigate two versions of the problem and present fully polynomial-time approximation schemes. 
In \cite{mia1}, the problem of scheduling deteriorating jobs with given release dates on a single batching machine is considered with the objective to minimize maximum lateness. After proving $NP$-hardness, a 2-approximation algorithm is given for the case of negative due dates of all jobs. It is also proven that the ED rule delivers an optimal solution for the case when all jobs have agreeable release dates, due dates and deterioration rates.  Miao et al. \cite{mia2} consider parallel and single machine problems with step-deteriorating jobs, given release dates and minimizing the makespan. For the single machine problem, NP-hardness in the strong sense has been proven. Finally, we mention that a very recent survey about deterioration (and learning) effects in scheduling has been given by Pei et al. \cite{pei}.

\section{Preliminaries and basic properties}

We use the existing notions, terminology and known properties that we briefly 
overview in this section. We refer the 
reader to, e.g., \cite{1div}, for more details and illustrative examples. 

A commonly used method for scheduling jobs with due dates was proposed by 
Jackson \cite{J} for the model without release times. Later it was 
extended by Schrage \cite{Sc} for non-simultaneously released jobs. 
The algorithm, iteratively,  among all currently released jobs, schedules
a most {\em urgent} job, one with the minimum due date; it updates
the current time $t$ and repeats the same operation for the updated current time: 
initially,  $t:=\min_j r_j$. Then, iteratively, $t$ is set to the minimum among the
release time of a yet unscheduled job and the completion time of the last assigned
job to the machine. 

We will refer to a schedule
generated by the heuristic as an {\em ED-schedule}, and denote the ED-schedule
constructed for the original problem instance (with the original job processing 
times) by $\sigma$. We may easily observe that 
in an ED-schedule $S$, there will arise an idle-time interval 
or a {\it gap} if the release time of a yet unscheduled job is strictly larger than 
the completion time of the job assigned as the last one to the machine. 

One can look at the schedule $S$ as a sequence of {\it blocks}, 
a consecutive sequence of jobs such that for each two neighboring jobs $i$ and 
$j$, job $j$ starts directly at the completion time of job $i$ behind time $r_j$.
If job $j$ starts at time $r_j$ then job $j$ starts a new block; in this case, 
we will say that a 0-length gap between the two jobs occurs.

We will distinguish a ``most critical'' part in a schedule $S$, one containing a 
job that realizes the value of the objective function. Let us call such a job an 
{\em overflow job}; i.e., if $o$ is an overflow job in the schedule $S$, 
then $$L_o(S)=\max_j L_j(S).$$ In case there are two or more consecutively 
scheduled jobs realizing the maximum objective function value, the last one 
is set to be the overflow job. 

\smallskip

A critical segment in the schedule $S$ containing an overflow job $o$ is called a 
{\it kernel} in that schedule:  it is a longest consecutive sequence of jobs ending
with job $o$ such that no job from this sequence has a due date larger than 
$d_o$. By definition, there may exist no gap in a kernel, and a kernel
is contained within some block. Since there may exist more than one overflow job, 
there will be the same amount of kernels in the schedule $S$. We use 
$$r(K)=\min_{i\in K}\{r_i\}$$ for the minimum release time of a job of kernel $K$, 
and abusing the notation, we will also use $K$ for the set of jobs from the kernel $K$.

Suppose the first job of kernel $K$ starts behind time $r(K)$ in schedule $S$. 
Then note that the completion time of the overflow job $o\in K$ can be decreased if 
some job $j\in K$ gets scheduled earlier than the earliest job of kernel $K$ was
scheduled in $S$. 
We call job $e$ scheduled immediately before the earliest scheduled job of kernel $K$ 
the {\em delaying emerging} job of kernel $K$ in schedule $S$. Note that job $e$
causes a forced right-shift for the jobs of the kernel $K$ in the schedule 
$S$, and that $d_e > d_j$. We will use $E(S)$ for the set of the
delaying emerging jobs in the schedule $S$.

\begin{Property} (Proposition 1 in \cite{1div})\label{opt}
If ED-schedule $S$  contains a kernel with no delaying emerging job, then it is optimal. 
\end{Property}

Assume a kernel $K\in S$ possesses the delaying emerging job $e\in E(S)$, and 
let $$\Delta(K) = f_e(S)-r(K)$$ be the forced right-shift or the {\em delay} for kernel 
$K$ in the schedule $S$. Furthermore, let $$\Delta_{\min}(S)=\min_{K\in S} \Delta(K).$$ 
The following property states previously known facts (e.g., see \cite{1div}).

\begin{Property}\label{difference}  
Given $K\in S$ with the delaying emerging job $e$ and the overflow 
job $o$,  we have $$L_o(S)-L^{\opt}\le \Delta_{\min}(S)\le \Delta(K) < p_{\max} .$$ 

\end{Property}

\comment
Observation. distinguish two kernels if they belong to different blocks: if 
not (there occurs pushing in one block) then by decreasing in the first kernel, 
in the second and following one will decrease, if gap occurs (cannot be decreased
more), then ok we have already two block...

---------------------
\endcomment

\subsection{Regularizing Kernels}

Based on Property \ref{opt}, it remains to consider the case, where
every kernel from the ED-schedule $\sigma$ possesses the delaying emerging job. 
Let $e$ be the delaying emerging job for kernel $K$ with the overflow job $o$ 
in an ED-schedule $S$. And let us consider an auxiliary partial ED-schedule 
$\bar K$ constructed solely for the jobs of the kernel $K$ by ED-heuristic. In 
this partial schedule, the earliest included job $i$ of the kernel $K$ starts at 
its release time $r_i=r(K)$. We shall refer to the kernel $K$ as {\it regular} if 
the overflow job in the partial schedule $\bar K$ has the due date $d_o$ (we 
note that the kernel in the schedule $\bar K$ will not be formed by all jobs 
of the kernel $K$ in case a gap in between them in that schedule occurs). 

As it is easy to see, if a kernel $K\in S$ is {\em irregular} (i.e., it is not 
regular), then the processing order of at least one pair of jobs $(i,j)$ in the 
schedule $S$ is reverted in the schedule 
$\bar K$. In other words, there occurs an {\em anticipating job} $j$, a job 
from kernel $K$ which ordinal number in the schedule $\bar K$ is less, compared 
to that in the schedule $S$. Since job $j$ is rescheduled before job $i$ in 
schedule $\bar K$, $r_j\le r_i$ (and $r_i > r(K)$), and also
$d_j > d_i$ as ED-heuristic included $i$ before $j$ in the schedule $S$. 

It follows that a former kernel job $j$ becomes the anticipated delaying 
emerging job for a newly arisen kernel in the schedule $\bar K$ in case
$K$ is irregular. The set of irregular kernels is similarly determined
in the schedule $\bar K$, and a similar construction is carried out for each of
them. This results in the omission of the newly arisen anticipated delaying 
emerging jobs. A recursive $O(|K|^2\log |K|)$ time procedure {\em fully 
decomposes (collapses)} the kernel $K$ into a sequence of so-called {\em 
substructure components} of that kernel with no irregular kernel.  Note that
the union of the set of jobs from these components is the set of jobs from the 
kernel $K$ minus the set of all the omitted anticipated delaying emerging jobs. 
We refer the reader to Section 4 of \cite{1div} for a detailed description of 
the decomposition procedure and examples illustrating it.

\begin{Property}\label{ED}
Suppose $K$ is a kernel in a schedule $S$. Then by increasing the 
processing time of the delaying emerging job $e$ of that kernel by $\tau$ time 
units, $0\le\tau\le\Delta(K)$, the maximum job lateness of a job in kernel $K$ 
will be increased by at most $\tau$ time units. If $K$ is a regular kernel, 
then by decreasing the processing time of job $e$ by $\tau$ time units, either 
the maximum job lateness of a job in kernel $K$ will decrease by $\tau$ time units 
or/and it will become a valid lower bound on the optimum job lateness. 
\end{Property}
Proof. First, we note that since a kernel may contain no emerging job, the jobs 
of any kernel $K$ are processed in an ED-order in schedule $S$, except that 
these jobs may be processed in a {\it weak} ED-order if there is a job $j$ with 
$d_j=d_o$ in the kernel $K$, where $o$ is the overflow job in that kernel: since 
$j$ is  a non-emerging job, it will not ``split'' the kernel $K$; it may start the 
kernel or may be included somewhere in between some more urgent jobs of the kernel
and hence the processing order will not be strictly ED. 
Suppose there exists such a job $j$ and the processing time of the delaying emerging 
job $e$ is increased. Then the time moment when job $j$ will be considered is 
respectively increased. As a result, a more urgent job of the kernel $K$, which 
was not yet released by the time moment $s_j(S)$ may get released and hence be
included by ED-heuristic ahead job $j$. Applying the same reasoning to all jobs $j$ 
with $d_j=d_o$, we obtain that the last scheduled job of the kernel $K$ with the 
due date $d_o$ will be the overflow job in the modified schedule and will
start $\tau$ time units later than $s_o(S)$. Similar reasoning holds for the first 
claim in case there exists no job $j$ with $d_j=d_o$, since an ED processing 
order of the jobs of the kernel $K$ will be kept in the modified schedule. 
This shows the first claim in the property.

For the second claim, assume that the processing time  of job $e$ is decreased by  
$\Delta$ time units; i.e., the first job of the kernel $K$ starts at its release 
time. Note that in the resultant complete ED-schedule $S'$, all jobs of the kernel 
$K$ will be scheduled as they were scheduled in the partial schedule $\bar K$. 
Furthermore, since the kernel $K$ is regular,  schedules $S$, $\bar K$ and 
$S'$ possess the same overflow job $o$. Then the second 
claim immediately follows if there arises no gap in between the jobs of the kernel 
$K$ in the schedule $\bar K$ ($S'$). Suppose there arises a gap, i.e., the kernel $K$
decomposes into substructure components. Then there is a single kernel $K'$ in schedule 
$\bar K$ that belongs to its last substructure component (see Lemma 4 in \cite{1div}). 
Since kernel $K$ is regular, there may exist no anticipated job and hence no emerging 
job for the kernel $K'$. Then the lateness of job $o$ is a lower bound on the
minimum jobs lateness (see Lemma 3 in \cite{1div}). Now for any $\tau<\Delta$ 
the same reasoning can obviously be applied, which completes the proof.\eop 

\smallskip

{\bf Kernel regularization procedure.} 
Now we describe a procedure that regularizes all irregular kernels an ED-schedule $S$.
First, it carries out the decomposition procedure for every irregular kernel $K\in S$. 
Let $D$ be the set of the delaying emerging jobs omitted during the
decomposition of all irregular kernels, and let $K^*$ be the  partial schedule, 
formed by the substructure components in the full decomposition of the kernel $K$. 
Note that schedule $K^*$ initiates at time $r(K)$, and that 
$\cup_{K\in S} K^* \cap D=\emptyset$ and $\cup_{K\in S} K^* \cup D = J.$  
We compose a complete schedule $(S)^*$ from the schedule $S$ in two passes. Pass 1
first  copies the schedule $S$, except its parts containing irregular kernels 
and the corresponding delaying emerging job, into an auxiliary partial schedule. 
The latter auxiliary schedule is completed by the partial schedules $K^*$ for each 
irregular kernel $K\in S$. Note that the resultant partial schedule of pass 1 is
feasible, i.e., there will occur no overlapping with the jobs from the initial 
auxiliary partial schedule in that schedule.
Furthermore, there will arise a gap before every of these irregular kernels in that
partial schedule. This schedule is augmented with the omitted jobs from the set $D$
at  pass 2. These jobs are included by the ED-heuristic, as follows. Whenever 
there is a yet unscheduled already released job from the set $D$, a most urgent one
is included at the beginning of the earliest available gap and the following 
jobs are correspondingly shifted to the right (in their actual processing order).  

\begin{Theorem}\label{regular}
Every kernel possessing the delaying emerging job in the schedule 
$(S)^*$ is regular, and this schedule can be constructed in 
$$\sum_{K\in S} O(|K|^2\log |K|)+O(n\log n)$$
time.
\end{Theorem} 
Proof. As to the time complexity, the decomposition of each irregular
kernel $K$ takes $O(|K|^2\log |K|)$ time (see Theorem 1 in \cite{1div}), and 
the insertion of the jobs from the set $D$ will take an additional time of 
$O(n\log n)$. As to the first claim in the theorem, observe that, all
kernels in the partial schedule of pass 1 are regular. Hence it will
suffice to show that the insertion of no delaying emerging job $e\in D$ 
will cause the rise of an irregular kernel. Suppose thus, once included, job $e$ 
yielded a forced right-shift for the subsequently scheduled jobs, that, 
in turn,  provoked the rise of a new kernel $K'$ in the schedule $(S)^*$.  
We need to show that the kernel $K'$ is regular. By the way of contradiction, 
suppose it is irregular, and let $j$ be the first anticipated job of that kernel. 
Let $i$ be a job preceding a job $j$ in the schedule $S$, such that job $j$ was 
included ahead job $i$ in the schedule $(S)^*$. We now observe that, because of the 
right-shift caused by job $e$,  job $j$ cannot start in the schedule $(S)^*$ before the 
time moment $s_i(S)$. Then both jobs should have been ready by the time moment, when 
job $j$ was included in the schedule $(S)^*$. Then the ED-heuristic could not include 
job $j$ before job $i$ since $d_j>d_i$, a contradiction. Thus, there may exist no 
anticipated job for any newly arisen kernel and hence it is regular.\eop

\section{The algorithms} 

\subsection{The first polynomial-time algorithm}

By Property \ref{opt}, if there is a kernel in the schedule $(\sigma)^*$ with no 
delaying emerging job, then it is optimal. Hence, assume from here on that 
every regular kernel in schedule $(\sigma)^*$ possesses the delaying emerging job. 
Then by reducing the processing time of the delaying emerging job of each regular 
kernel of this schedule, the current maximum job lateness can be reduced:
Algorithm 1 compresses the delaying emerging job of each kernel 
$K\in(\sigma)^*$ by the amount $\Delta_{\min}(\sigma)$ and shifts the jobs 
succeeding each compressed delaying emerging job to the left correspondingly.
We denote the resultant schedule by  $\sigma_-\Delta_{\min}$ (for notational 
simplicity we omit argument $\sigma$ in $\Delta_{\min}$). Note that in  schedule
$\sigma_-\Delta_{\min}$, the former overflow job of every kernel from schedule 
$(\sigma)^*$ has the same lateness.

\comment
Not necessarily the  schedule $\sigma_-\Delta_{\min}$ satisfies the feasibility 
condition (\ref{1}). Moreover,   
a new overflow job may occur in the schedule $\sigma_-\Delta_{\min}$. Then if
the maximum job lateness in the schedule $\sigma_-\Delta_{\min}$ is realized by 
an overflow job from the schedule $(\sigma)^*$, i.e., there arises no new kernel 
in the former schedule, and, in addition, schedule $\sigma_-\Delta_{\min}$ 
does not violate the feasibility restriction (\ref{1}), it is optimal:
\endcomment

\begin{Theorem}\label{aux 1}
Algorithm 1 constructs the schedule $\sigma_-\Delta_{\min}$ in time
$$\sum_{K\in (\sigma)^*} O(|K|^2\log |K|)+O(n\log n),$$
and it is an optimal feasible solution if \\[1ex] 
(i) $$\sum_{e\in E((\sigma)^*)} c_e \Delta_{\min}((\sigma)^*) \le U,$$ and  \\[1ex]
(ii) it contains the same set of kernels as the schedule $(\sigma)^*$. 
\end{Theorem} 
Proof. As to the construction cost, the construction of the original schedule 
$\sigma$ takes $O(n\log n)$ time, and the regularization of each kernel in that 
schedule and the creation of the schedule $(\sigma)^*$ has an additional cost of 
$$\sum_{K\in (\sigma)^*} O(|K|^2\log |K|)+O(n\log n)$$ (Theorem \ref{regular}). 
The required parameters of each of these regular kernels from the schedule 
$(\sigma)^*$ can 
be obtained in $O(n)$ time, and the left-shift of the corresponding schedule 
portions after the compression of each delaying emerging job from the set 
$E((\sigma)^*)$ can also be accomplished in $O(n)$ time since the processing 
order of the jobs need not be changed (see again Property \ref{ED}). 

As to the optimality, suppose the feasibility condition (i) is satisfied. 
Since the condition (ii) is also satisfied, there may exist no anticipated 
job for a kernel in schedule $\sigma_-\Delta_{\min}$ and hence the processing 
order of the jobs of every kernel is optimal (see the proof of Property \ref{ED}).
Furthermore, let $K$ be a kernel in schedule $\sigma_-\Delta_{\min}$
with $$\Delta(K)=\Delta_{\min}(\sigma_-\Delta_{\min}).$$ The first job of 
every kernel $K$ with this property starts at its release time in the schedule 
$\sigma_-\Delta_{\min}$. Then the maximum job lateness is minimized in that 
schedule since the processing order of the jobs of every kernel agrees
with an optimal sequence.\eop

\comment
Now, it is a known fact that, 
for a given set of jobs, among all partial schedules that start at the same time, 
one in which the jobs from this set appear in ED-order minimizes the maximum job
lateness. Indeed, using an interchange argument,  we can easily verify that there 
is no benefit in interchanging  any two jobs from that partial schedule. Now the
theorem follows since this reasoning can be applied to every kernel from the 
schedule $\sigma_-\Delta_{\min}$.
\endcomment

\subsection{The pseudo-polynomial-time algorithm}

It remains to study the case when the schedule $\sigma_-\Delta_{\min}$ does not
satisfy the conditions from Theorem \ref{aux 1}. Assume that a new 
kernel occurs in the schedule $\sigma_-\Delta_{\min}$ (see condition (ii) 
in Theorem \ref{aux 1}). Algorithm 2 regularizes every such kernel 
and then compresses the processing time of the corresponding delaying 
emerging jobs by just one time unit. It  repeats the same operations for each 
newly created schedule, as long as this schedule remains feasible 
(i.e., it satisfies condition (\ref{1})): \\[0.5ex]

\noindent
{\bf Algorithm 2.}  \\[1ex]

\noindent {\bf Step 0.}  \\[0.5ex]

 Create the initial ED-schedule $\sigma$; $S:=(\sigma)^*$; \\[0.5ex] 
 
\noindent {\bf Step 1.} \\[0.5ex]

IF all kernels in the schedule $S$ possess the delaying emerging job  \\[0.5ex]

\hskip.2cm  THEN create a modified auxiliary ED-schedule $S'$ from the schedule $S$
  by compressing 
       
      \hskip1.6 cm the processing time of the  delaying emerging job of each 
      kernel in schedule $S$  
      
      \hskip1.6 cm by  one time unit  and shifting the 
      following jobs correspondingly to the left; \\[0,5ex]
       
  \hskip1.2cm    IF in the schedule $S'$ the feasibility condition (\ref{1}) is not 
  violated   \\[0.5ex]
  
   \hskip1.8cm  THEN  $S:=(S')^*$ \{regularize schedule $S'$\};  repeat Step 1 \\[0.5ex]
   
     \hskip1.8cm  ELSE return schedule $S$ \{$S$ is a feasible schedule\} \\[0.5ex]
     
\hskip.2cm ELSE return schedule $S$.

\bigskip

We will refer to a schedule $S$ as {\em well-balanced} if by dis-compressing 
the delaying emerging job of all kernels of that schedule by $\tau$ time units, the 
maximum job lateness will increase by the same amount. 
Note that the schedule $S$ created by Algorithm 2 has a nice property that the lateness 
of the overflow jobs of all kernels in it is the same.  As a consequence, this schedule 
is well-balanced (see Property \ref{ED}): 

\begin{Property}\label{well-balanced}
The schedule $S$ of every iteration $h$ in Algorithm 2 is well-balanced.
In particular, the maximum job lateness in that schedule is simultaneously 
attained in all kernels from the set ${\cal K}^h$. 
\end{Property}

Let us also observe that every  kernel of iteration $h-1$ will remain 
regular at iteration $h$ and the following iterations (again Property \ref{ED}, 
see also Theorem \ref{regular}). At the  same time,  since the lateness of the 
overflow jobs of iteration $h-1$ are decreased, a new kernel may arise at 
iteration $h$ in the schedule $S'$. If a newly arisen kernel in schedule $S'$ of 
iteration $h$ is irregular, it is regularized at iteration $h$ in the newly created 
schedule $S=(S')^*$. We will use  ${\cal K}^h$ for the set of all the  (regularized)
kernels by  iteration $h$; as we just observed, 
${\cal K}^h\subseteq {\cal K}^{h-1}$ holds for any $h>1$ in Algorithm 2.

\begin{Property}\label{one less}
The maximum job lateness in the schedule $S$ of each iteration $h$ in Algorithm 2 
(except possibly the last iteration) is exactly one less than that of iteration $h-1$.  
\end{Property} 
Proof. By the construction of step 1, at every iteration $h>1$, the maximum job 
lateness in the schedule $S'$ is one less than that in the regularized schedule $S$ 
of iteration $h-1$. Then by Property \ref{ED}, it suffices to show that the maximum 
job lateness in the regularized schedule $S$ of iteration $h$ is the same as that 
in the schedule $S'$ of the same iteration. There are two possible cases. If  
there arises no new kernel in the schedule $S'$ of iteration $h$, then this
is obviously true. Otherwise, note that the overflow job of a newly arisen 
kernel will have the lateness equal to the maximum job lateness in the schedule 
of iteration $h-1$ minus one, and the maximum job lateness in the regularized
schedule $S$ must be the same as that in the schedule $S'$.\eop 

\begin{Theorem}\label{opt-1} 
Algorithm 2 creates a well-balanced feasible solution $S$ for problem 
$1|r_j, compressible(U)|L_{\max}$ in less than $p_{\max}$ iterations. 
If it does not create an infeasible solution $S'$ (the first ELSE condition
is never entered at Step 1), then $S$ is also optimal. Otherwise
(the first ELSE condition is executed at Step 1 before the algorithm halts), 
the solution $S$ is optimal for the generic problem $1|r_j|L_{\max}$ with the 
compressed set of job processing times defined by the compression vector 
$(x_1(S),\dots,x_n(S))$. 
\end{Theorem}
Proof. By Property \ref{one less}, the maximum job lateness in the schedule $S$ 
of each (non-terminal) iteration $h$ in Algorithm 2 is one less than 
that of iteration $h-1$. Then the number of iterations is bounded by $p_{\max}$ 
due to Property \ref{difference}. The schedule $S$ is well-balanced by Property 
\ref{well-balanced} and it satisfies feasibility condition (\ref{1}
by the construction of step 1. 

Suppose now at step 1 the first ``ELSE'' statement is never entered.   
Then the feasibility condition (\ref{1}) was never violated and the algorithm 
has stopped by entering the second ``ELSE'' statement; i.e., there is a kernel 
in  schedule $S$ possessing no delaying emerging job and this schedule $S$ is 
optimal  by Property \ref{opt}.

Suppose now that the algorithm halts by entering the first ``ELSE'' statement, 
i.e., the feasibility condition (\ref{1}) does not hold in the schedule $S'$
and the schedule $S$, regularized at the previous iteration,  is returned. We
have to show that $S$ is an optimal feasible schedule. Similarly as in the previous
case, $S$ is feasible by the construction of step 1. Furthermore, since all kernels
in that schedule are regular,  by changing the processing  order of the jobs
of any  kernel, the lateness of the corresponding overflow job cannot be 
decreased. Let us now consider the remaining non-kernel jobs of the schedule. Without
loss of generality, assume that schedule $S$ consists of a single block, since 
our reasoning can be applied to each individual block independently. Schedule $S$ 
consists of kernels and alternative sequences of jobs scheduled by ED-heuristic 
before the first kernel, in between two neighboring kernels, and after the last kernel. 
Since there is no gap within any of these sequences, changing the processing order 
of the jobs from any of these sequences will either leave the current lateness 
unaltered or will increase it (due to a possible gap that may occur because of the
order change). Interchanging the jobs from different sequences may again cause
new gaps. As a result, if the lateness of some overflow job from the schedule 
$S$ decreases, that of some other overflow job will increase, see Properties 
\ref{ED}, \ref{well-balanced} and \ref{one less}. It follows that schedule $S$ 
minimizes maximum job lateness for the problem instance of the generic 
problem $1|r_j|L_{\max}$ with the compressed set of job processing times
determined by Algorithm 2.\eop

\subsection{The second polynomial-time algorithm}

In the previous sub-section, we constructed a feasible schedule which is optimal
for the instance of the generic problem $1|r_j|L_{\max}$ with the compressed set 
of job processing times determined by Algorithm 2. Since the compression at each 
iteration is carried out by one unit of time, the total number of iterations
depends on the maximum job processing time and hence the algorithm runs in
pseudo-polynomial time.  At first glance, $O(p_{\max})$  seems to be a natural 
bound on the total number of iterations since the lateness of the overflow job 
in a newly arisen kernel may differ just by one time unit from that of an 
earlier arisen kernel. Nevertheless,  maintaining schedules 
with regular kernels is helpful and will permit us to obtain a schedule with 
the same property as the one delivered by Algorithm 2 in $O(n)$ iterations. 

This section's polynomial-time Algorithm 3 combines some features of Algorithms 1 
and 2. At stage 1, it outputs a well-balanced schedule, which is optimal
for the instance of the generic problem $1|r_j|L_{\max}$ with the 
obtained compressed set of job processing times. If this solution is
not feasible for problem $1|r_j, compressible(U)|L_{\max}$, then at stage 2, 
the solution of stage 1 is converted to 
a well-balanced feasible schedule. Like Algorithm 1 (and unlike Algorithm 2),
Algorithm 3 makes jumps of $\Delta_{\min}(\sigma_{h-1})$ while compressing the 
processing times of the delaying emerging jobs at iteration $h$ at stage 1, 
which guarantees polynomial-time performance. As Algorithm 2, Algorithm 3 
maintains the schedule $\sigma_h$ of each iteration $h$ regular and well-balanced. 
Since, unlike Algorithm 2, the ``compression jumps'' 
in Algorithm 3 are carried out by more than one units time, a special 
care is to be taken to keep the schedule of every iteration well-balanced
(note that, unlike Algorithm 2, the maximum job lateness in a non-regularized and 
the corresponding regularized schedules will not necessarily be equal).

\subsubsection{Stage 1}

At the initial iteration 1, stage 1 invokes Algorithm 1, and the schedule 
$\sigma_1=\sigma_-\Delta_{\min}$ is obtained. For convenience, 
let $\sigma_0=\sigma$, and let  
$$\Lambda^1=L(\sigma_0)-\Delta_{\min}(\sigma_0) = L_o(\sigma_1)$$ 
be the lateness of an overflow job $o\in \sigma_0$ in the schedule 
$\sigma_1$ (recall that
schedule $\sigma_-\Delta_{\min}$  is obtained from a regularized 
schedule $(\sigma)^*$ by compressing the processing time of each delaying
emerging job of that schedule by the amount $\Delta_{\min}(\sigma)$, so that 
the lateness of the overflow jobs of the schedule $\sigma=\sigma_0$ is decreased 
by $\Delta_{\min}(\sigma)$ time units in the schedule $\sigma_1$, see Property \ref{ED}). 

If $\Lambda^1 = L(\sigma_1)$, i.e.,  
in the schedule $\sigma_1$ no new kernel arises (condition (ii) in 
Theorem \ref{aux 1}),  stage 1 outputs this well-balanced (but not 
necessarily feasible) schedule. Note that $\Lambda^1 > L(\sigma_1)$ is not 
possible. If now $\Lambda^1 < L(\sigma_1)$ 
(at iteration $1$, a new overflow job arises in the schedule $\sigma_{1}$), 
stage 1 proceeds with the next iteration 2. Again, first, schedule $\sigma_1$ is 
regularized, and then the processing time of every delaying emerging job in the 
regularized schedule $(\sigma_{1})^*$ is compressed by the magnitude 
$\Delta_{\min}((\sigma_{1})^*)$ and the jobs succeeding each compressed  
delaying emerging job are shifted to the left. This 
results in an auxiliary schedule ${\sigma_1}_-\Delta_{\min}$, 
in which each newly arisen kernel $K\in {\cal K}^2$ is delayed by 
$$\Delta(K)-\Delta_{\min}(\sigma_{1})$$ time units; in particular, a kernel 
$K$ with $$\Delta(K)=\Delta_{\min}(\sigma_{1})$$ starts at time $r(K)$ 
without any delay. Below we describe how delaying emerging jobs are 
iteratively compressed and dis-compressed.

\smallskip

In general, for a given iteration $h\ge 1$, let  
$$\Lambda^h = L(\sigma_{h-1})-\Delta_{\min}(\sigma_{h-1})=L_o(\sigma_h),$$ 
where $o$ is an overflow job in the schedule $\sigma_{h-1}$. Stage 1 
outputs the current schedule $\sigma_h$ if $\Lambda^h = L(\sigma_h)$, i.e., 
no new kernel in the schedule $\sigma_{h}$ arises. 

Suppose now $\Lambda^h < L(\sigma_h)$, i.e., a new kernel in the schedule 
$\sigma_{h}$ arises (again, $\Lambda^h > L(\sigma_h)$ is not possible). 
Then stage 1 proceeds with iteration $h+1$ by 
first creating an auxiliary schedule ${(\sigma_h)^*}_-\Delta_{\min}$. 
The following two cases are distinguished: (1) 
$$L ({(\sigma_h)^*}_-\Delta_{\min}) >\Lambda^h,$$ i.e., the 
lateness of an overflow job in the schedule 
${(\sigma_h)^*}_-\Delta_{\min}$ is more than the 
maximum job lateness at iterations $1,\dots,h$. And (2) 
$$L ({(\sigma_h)^*}_-\Delta_{\min}) < \Lambda^h,$$ i.e.,
the lateness of an overflow job in the schedule 
${(\sigma_h)^*}_-\Delta_{\min}$ is less than the 
maximum job lateness of iterations $1,\dots,h$.

\begin{Property}\label{dis-a}
Suppose $$L ({(\sigma_h)^*}_-\Delta_{\min}) >\Lambda^h.$$ Then 
$$L^{\opt} \ge L ({(\sigma_h)^*}_-\Delta_{\min}).$$ Hence, 
the delaying emerging job of every kernel in ${\cal K}_{h}$ can be 
dis-compressed by $$L ({(\sigma_h)^*}_-\Delta_{\min})- \Lambda^h$$ time 
units in the schedule ${(\sigma_h)^*}_-\Delta_{\min}$.  
\end{Property}
Proof. Note that the schedule $\sigma_h$ is regularized at iteration $h+1$ 
and hence all kernels in ${\cal K}_{h+1}$ are regular. In particular, 
there is a regular kernel $K$ with $\Delta(K)=\Delta_{\min}(\sigma_{h})$ in the
schedule ${(\sigma_h)^*}_-\Delta_{\min}$. This kernel possesses no delaying emerging 
job and hence 
$$L^{\opt} \ge L_o({(\sigma_h)^*}_-\Delta_{\min})=L({(\sigma_h)^*}_-\Delta_{\min}),$$
where $o$ is the overflow job in the kernel $K$. The second claim in the property 
now obviously follows.\eop  


\begin{Property}\label{dis-b}
$L^{\opt} \ge \Lambda^h$ holds. Hence, if 
$$L ({(\sigma_h)^*}_-\Delta_{\min}) < \Lambda^h,$$ then the delaying emerging 
job of every (new) kernel in ${\cal K}^{h+1}\setminus {\cal K}^{h}$ can be 
dis-compressed by $$\Lambda^h - L ({(\sigma_h)^*}_-\Delta_{\min})$$ time 
units in the schedule ${(\sigma_h)^*}_-\Delta_{\min}$.  
\end{Property}
Proof. Similar to that of Property \ref{dis-a}.\eop

Stage 1  dis-compresses delaying emerging jobs according to Properties \ref{dis-a}
and \ref{dis-b}, shifting to the right the jobs succeeding each dis-compressed delaying 
emerging job correspondingly. This results in the schedule $\sigma_{h+1}$. In case (1)  
$L(\sigma_{h+1})= L_o ({(\sigma_h)^*}_-\Delta_{\min})$, as indicated earlier, and in case 
(2) $L(\sigma_{h+1})= \Lambda^h$ (the current maximum job lateness is kept 
unchanged). These operations are carried out as long as in schedule $\sigma_{h+1}$
a new kernel/overflow job arises. Otherwise, stage 1 outputs a well-balanced 
schedule $\sigma_{h+1}$. \\[0.5ex]

\noindent  
{\bf Algorithm 3, Stage 1.}  \\[1ex]

\noindent {\bf Step 0.}  \\[0.5ex]

 Create the initial ED-schedule $\sigma$ and regularize it; 
 
 Create the schedule $\sigma_1=\sigma_-\Delta_{\min}$; $h:=1$; \\[0.5ex] 
 
\noindent {\bf Step 1.} \\[0.5ex]

IF $\Lambda^h = L(\sigma_h)$ \{no new kernel arises in the schedule $\sigma_h$\} \\[0.5ex]

\hskip.2cm  THEN  RETURN schedule $\sigma_h$ \\[0.5ex]

\hskip.2cm ELSE \{$\Lambda^h < L(\sigma_h)$\} \\[0.5ex] 
  
\hskip1.2cm IF $L ({(\sigma_h)^*}_-\Delta_{\min}) >\Lambda^h$
\{ the current maximum job lateness is surpassed \} \\[0.5ex]

\hskip1.8cm  THEN determine the schedule $\sigma_{h+1}$ by dis-compressing in the
schedule  

\hskip2cm ${(\sigma_h)^*}_-\Delta_{\min}$ the processing time of the  
delaying emerging job  

\hskip2cm of every kernel in ${\cal K}_{h}$  by 
$L ({(\sigma_h)^*}_-\Delta_{\min})- \Lambda^h$ time units \\[0.5ex]

\hskip1.8cm ELSE determine the schedule $\sigma_{h+1}$ by dis-compressing in the 
schedule  

\hskip2cm ${(\sigma_h)^*}_-\Delta_{\min}$ the processing time of the delaying 
emerging job of every

\hskip2cm  (newly arisen) kernel in ${\cal K}^{h+1}\setminus {\cal K}^{h}$  by 
$\Lambda^h - L ({(\sigma_h)^*}_-\Delta_{\min})$ time units; \\[0.5ex]   

\hskip.2cm $h:=h+1$; REPEAT Step 1.        
      
\bigskip

\comment
\begin{Property}\label{lb 2}
At any iteration $h$, $\{ L_o(\sigma_{h-1})-
The lateness of any overflow job from  one of the kernels from set 
${\cal K}_h$ is a lower bound on the optimum schedule lateness. In particular, 
$$\min_h \{ L_o(\sigma_{h-1}) | o\in K, K\in {\cal K}_h\}$$ is a valid lower
bound.
\end{Property}
Proof.  By the construction of the algorithm, the completion time of the
overflow job $o$ of kernel $K\in {\cal K}_h$ is decreased by 
$\Delta(K)-\Delta_{\min}(\sigma_{h-1})$ time units so that the first job
of that kernel starts at time $r(K)$. 
\endcomment

\begin{Theorem}\label{alg 3}
Stage 1 delivers a well-balanced schedule $\sigma^{h}$ possessing a kernel 
without the delaying emerging job in $O(n)$ iterations in $O(n^2\log n)$ time. 
The schedule $\sigma^{h}$ is optimal for the corresponding instance of the generic 
problem  $1|r_j|L_{\max}$ with the compressed set of job processing times 
defined by the compression vector $(x_1(\sigma),\dots,x_n(\sigma))$, and 
it is also optimal for the original instance of problem 
$1|r_j, compressible(U)|L_{\max}$, if it does not violate the feasibility 
condition (\ref{1}). 
\end{Theorem}
Proof. Initially at step  0, the cost of invoking Algorithm 1 is as in 
Theorem \ref{aux 1}. Iteratively, in iteration $h>1$ at step 1, the regularization 
of each kernel $K\in {\cal K}^{h-1}$ from the schedule $\sigma_{h-1}$ has the cost of 
$O(|K|^2\log |K|)$ (Theorem \ref{regular}), and similarly as in Algorithm 1, the 
required parameters of each of these regular kernels can be obtained in $O(n)$ 
time. After each compression (dis-compression, 
respectively) of the delaying emerging jobs, the left-shift (right-shift, 
respectively) of the jobs of the corresponding kernels needs $O(n)$ time. 

Assume, for now, that the total number of different kernels 
that may occur during the execution of stage 1 is bounded from above by $O(n)$.
Then the above operations at step 1 have to be repeated a number of times
bounded by $O(n)$ and the bound $O(n^2\log n)$ follows.  Below 
we show that in less than $n$ iterations the algorithm arrives at an
iteration ${h}$ such that there occurs no new overflow job at that iteration 
(equivalently, ${\cal K}^{h}\setminus {\cal K}^{h-1} = \emptyset$, i.e.,
there occurs no new kernel, condition (ii) from Theorem \ref{aux 1}). 

Thus we wish to show that in no more than $n$ iterations, stage 1 arrives at an 
iteration ${h}$ with ${\cal K}^{h}\setminus {\cal K}^{h-1} = \emptyset$, as we 
claimed above. Consider an arbitrary iteration $g$ with a newly arisen kernel 
$K^g \in {\cal K}^g$, $K^g \not\in {\cal K}^{g-1}$. There may clearly occur 
no more than $n$ such iterations if ${K}^g\cap {K}^f=\emptyset$, for any 
$K^f \in {\cal K}^f$, $f=1,\dots,g-1$. If now ${K}^g\cap {K}^f\ne\emptyset$, 
then the two kernels $K^f$ and $K^g$ must be  
identical. Indeed, let, first, $j$ be any job scheduled after the kernel $K^f$. 
Note that job $j$ could not have been shifted to the right by more than any job 
of kernel $K^f$. Then job $j$ cannot become 
a part of kernel $K^g$ since it did not form a part of kernel $K^f$. Suppose that 
$j$ is a job scheduled before the delaying emerging job $e$ of kernel
$K^f$ in the schedule $\sigma_{f-1}$. If job $j$ belongs to kernel $K^g$ then
job $e$ must also belong to that kernel, but no  kernel may contain a 
(compressed) delaying emerging job. The kernel $K^g$ cannot form a 
sub-sequence/subset of the kernel $K^f$ since both kernels are regular. 
It follows that the two kernels coincide (as sequences of jobs), and  they
have the same overflow job $o$. But the lateness of job $o$ is optimal
since this kernel is regular and its first job starts at its release time. We 
showed that $g=h$ and the schedule $\sigma_h$ is optimal for the instance of problem 
$1|r_j|L_{\max}$ with the compressed set of job processing times. Note also that
the schedule $\sigma_h$ is well-balanced since, by the construction, the lateness 
of the overflow jobs of all the  kernels arisen at iterations $1,\dots,h$ is the 
same. The last claim in the theorem obviously follows. \eop

\comment
Now we turn to the second adjustment stage that 
dis-compresses (increases) the processing time of some compressed delaying 
emerging jobs. This results in the reduction of the total cost without increasing  
the current maximum lateness. We need a few additional notations to describe
the second stage. 

We will refer to a kernel  $K\in {\cal K}_h$ with 
$$\Delta(K)=\Delta_{\min}(\sigma_{h-1})$$ as a {\it critical kernel} in 
the schedule $\sigma_{h-1}$. Furthermore, we let $$L^+(h)=L_o(\sigma_{h-1}) \quad \mbox{and} \quad 
L^-(h)=L_o(\sigma_{h})$$ be the lateness of the overflow job $o$ of a critical kernel 
$K$ in the schedules $\sigma_{h-1}$ and $\sigma_h$, respectively. We may observe that, 
by the construction, $$L^+(h)-L^-(h) = x_e(\sigma_h),$$ where $e$ is the corresponding
delaying emerging job. Stage 1 will stop at iteration $h$ if 
$$L^-(h)=L_{\max}(\sigma_h).$$ Otherwise, $L^-(h) < L_{\max}(\sigma_h)$, and stage
1 will proceed with iteration $h+1$ with a new set of kernels ${\cal K}_{h+1}$.  

Now, we observe 
that the processing time of some former delaying emerging jobs can be increased 
without increasing $L_{\max}(\sigma_{h^*})$. In particular, it is not difficult 
to see that the current processing time of the delaying emerging job of every 
former kernel $K\in \cup \{{\cal K}_h | h=1,\dots,h^* \}$ with 
$$L_{\max}(K) < L_{\max}(\sigma_{h^*})$$ can be increased by the amount 
$$L_{\max}(\sigma_{h^*})-L_{\max}(K)$$ without increasing $L_{\max}(\sigma_{h^*})$.

The adjustment stage 2 carries out this update of the processing times of the
corresponding delaying emerging jobs in $O(n)$ time (stage 1 keeps the track of the 
data of the kernels from the set $\cup \{{\cal K}_h | h=1,\dots,h^*\}$ so that the update
is accomplished in $O(n)$ time at stage 2). Thus, the worst-case time complexity 
of Algorithm 3 is that of stage 1 from Theorem \ref{alg 3}. We will refer to
the schedule delivered by stage 2 again by $\hat\sigma$. 

\endcomment

\subsubsection{Stage 2}

Stage 2 transforms the solution $\hat\sigma=\sigma^h$ of stage 1 into a feasible 
solution $\bar\sigma$ to problem $1|r_j, compressible(U)|L_{\max}$ by dis-compressing
the processing time of the compressed delaying emerging jobs. Observe that
$$\sum_{e\in {\cal E}}x_j(\hat\sigma) c_j - U$$ is the excess of 
the total compression cost in the solution $\hat\sigma$. Let $k=|{\cal K}^h|$ 
(recall that ${\cal K}^h$ is the set of all (regular) kernels formed at stage 1),
and let ${\cal E}$ be the set of the corresponding delaying emerging jobs. 
By increasing the processing time of the delaying emerging job $e\in {\cal E}$
of each kernel from the set ${\cal K}^h$ by 
$$\xi = \lceil (\sum_{e\in {\cal E}}x_j(\hat\sigma) c_j - U)/ k\rceil$$ time
units and shifting the jobs of that kernel correspondingly to the right, stage 2 
obtains the schedule $\bar\sigma$ in $O(n)$ time, in which the maximum job 
lateness is increased by $\xi$ time units compared to schedule  $\hat\sigma$.

\begin{Theorem} 
The schedule $\bar\sigma$ delivered by Algorithm 3 in $O(n^2\log n)$ time
is a feasible schedule for the problem $1|r_j, compressible(U)|L_{\max}$ and
satisfies the following equality:
\begin{equation}\label{xi}
L_{\max}(\bar\sigma)=L_{\max}(\hat\sigma)+\xi.  
\end{equation}
Furthermore, the schedule $\bar\sigma$ is well-balanced and optimal 
for the corresponding instance of the generic problem  $1|r_j|L_{\max}$ with the 
compressed set of job processing times in the schedule $\bar\sigma$. 
\end{Theorem}
Proof. The time complexity follows from Theorem \ref{alg 3} since stage 2 
clearly runs in $O(n)$ time. The feasibility of the schedule 
$\bar\sigma$ follows from the definition of the parameter $\xi$. Furthermore,
since every kernel $K\in \hat\sigma$ is regular, 
by increasing the processing time of the corresponding delaying emerging job by 
$\xi$ time units, the lateness of the corresponding overflow job will increase 
by the same amount (Property \ref{ED}). Then 
Equation (\ref{xi}) holds since the set of kernels in the schedule $\bar\sigma$ is 
the same as that in the schedule $\hat\sigma$. Moreover, the schedule $\bar\sigma$
is well-balanced since the schedule $\hat\sigma$ is well-balanced (Theorem \ref{alg 3}).
Then schedule $\bar\sigma$ is optimal for the instance of the problem  
$1|r_j|L_{\max}$ with the compressed set of processing times in that schedule,
similarly as in Theorems \ref{opt-1} and \ref{alg 3}.\eop

\comment

At every iteration, the lateness of every overflow
job is reduced to a lower bound on the optimum schedule lateness. Suppose $o$ is 
such an overflow job in schedule $\sigma_h$. By the 
construction, the corresponding delaying emerging job is compressed so that
----the lateness of job $o$ in the resultant schedule becomes a lower bound on
the optimum schedule lateness (--- see Property \ref{lb}). Moreover, this job cannot 
be rescheduled to a later time moment in any further generated schedule. Hence 
no further generated schedule may contain a kernel possessing the delaying emerging 
job in which job $o$ is an overflow job. Then the lateness of the overflow jobs of 
each iteration is strictly less than that of the overflow jobs of the previous
iteration(s). It follows that in at most $\Delta_{\min}(\sigma)$ iterations, there 
will arise a kernel which overflow job realizes a lower bound  on the optimum schedule
lateness and this kernel will have no delaying emerging job and the algorithm
will stop at that iteration. But since $\Delta_{\min}(\sigma) < p_{\max}$
(Property \ref{difference}) the bound $p_{\max}$ follows.   $n$ is another valid bound 
since, as we showed, the same job may become an overflow job in at most one
schedule such that job $o$ is an overflow job in it and
such that the corresponding kernel possesses the delaying emerging job.\eop 

At each iteration, the current maximum lateness
is decreased by at least one unit of time. Hence, the total number of times the
next created schedule $\sigma_h$ contains a kernel with the delaying emerging job
cannot exceed $\Delta_{\min}(\sigma)$ since $L_o(\sigma)-\Delta_{\min}(\sigma)$ 
(where $o$ is an overflow job in schedule $\sigma$) is a lower bound on the 
optimum schedule lateness and this bound is attained 
(Property \ref{lb}). i.e., in at most $\Delta_{\min}(\sigma)$
iterations there will arise a kernel which overflow job realizes this lower 
bound and this kernel will  

from schedule $\sigma$ will become a kernel 
But at least one  of these kernels possesses no delaying emerging job (by the
construction of schedule $\sigma_-\Delta_{\min}$). But  $\Delta_{\min}(\sigma)$
and the bound $p_{\max}$ follows.


 \endcomment

\section{Conclusion} 

We showed that by shortening the processing times of some specially determined 
(emerging)  jobs to ``requited amounts'' of time units, 
the maximum job lateness can be effectively reduced. 
Algorithm 3, in practice, provides the 
manufacturer with the information about the number of additional resources 
that are required to solve the problem  $1|r_j,compressed(U)|L_{\max}$ optimally. 
Such an approach is particularly useful in situations 
when the manufacturer is willing to allocate additional resources to shorten 
the processing of some late pending (emerging) jobs and provide in this way
a non-delay completion of recently arrived more urgent jobs. 
In the case when this amount of additional resources is not
available, stage 2 returns a  well-balanced schedule that 
is optimal for the corresponding instance of problem $1|r_j|L_{\max}$. 
Note that, a well-balanced schedule has an important optimal tightness property 
that, by dis-compressing the processing of any compressed delaying emerging job,
the maximum job lateness will be increased by the same amount. More generally, 
in a well-balanced schedule, by increasing the processing time of any kernel or 
non-kernel job preceding some kernel, the maximum job lateness will be increased by 
the same amount. Finally, 
the pseudo-polynomial time algorithm of Section 4.2 may be attractive 
for applications where $p_{\max}$ is a priory known non-large magnitude. 
For further research, it will be interesting to see how the proposed approach  
can be used for parallel machine environments.

\section{ Acknowledgment}

The research of the first author was supported by the DAAD grant 
57507438 for Research Stays for University Academics and Scientists.


\end{document}